\documentclass{amsart}
\usepackage{amsfonts}
\usepackage{amssymb}
\usepackage{amsmath}
\usepackage{amsrefs}
\usepackage{mathrsfs}

\usepackage{color}
\usepackage[dvipsnames]{xcolor}
\usepackage{tikz}
\usepackage{setspace}
\usepackage{multicol}
\usetikzlibrary{calc,intersections,decorations.pathreplacing}
\usepackage{ragged2e}
\usepackage[linktocpage=true]{hyperref}
\usepackage{pgfplots}

\usetikzlibrary{arrows.meta}
\usetikzlibrary{arrows}

\newtheorem{thm}{Theorem}[section]

\newtheorem{lemma}[thm]{Lemma}
\newtheorem{proposition}[thm]{Proposition}

\theoremstyle{definition}

\newtheorem{definition}[thm]{Definition}
\newtheorem{prop-def}[thm]{Proposition--Definition}
\newtheorem{claim}[thm]{Claim}

\newtheorem{term}[thm]{Terminology}

\theoremstyle{remark}
\newtheorem{rem}[thm]{Remark}
\newtheorem{question}[thm]{Question}
\newtheorem{example}[thm]{Example}

\usetikzlibrary{fit}

\hypersetup{ colorlinks=true, citecolor=cyan, linkcolor=blue, urlcolor=blue}

\begin{document}

\title{Deformation of quintic threefolds to the chordal variety}
\author{Adrian Zahariuc}
\email{azahariuc@math.ucdavis.edu}
\address{Harvard University, One Oxford Street, Cambridge MA, 02138, USA}
\curraddr{Department of Mathematics, UC Davis, One Shields Ave, Davis, CA 95616 }
\thanks{Partially supported by National Science Foundation grant DMS-1308244 \emph{Nonlinear Analysis on Sympletic, Complex Manifolds, General Relativity, and Graphs}.}

\begin{abstract}
We consider a family of quintic threefolds specializing to a certain reducible threefold. We describe the space of genus zero stable morphisms to the central fiber, as defined by J. Li. As an application of a straightforward extension, we prove the existence of rigid stable maps with smooth source of arbitrary genus and sufficiently high degree to very general quintics.
\end{abstract}

\maketitle

\bibliographystyle{amsplain}

\vspace*{6pt}\tableofcontents  
% A table of contents should normally not be included

\section{Introduction}
\label{sec:introduction}

\subsection{Algebraic curves on Calabi-Yau threefolds}

In the study of algebraic curves on Calabi-Yau manifolds, one of the obstacles to using naive geometric approaches is the fact that such curves can rarely be deformed holomorphically. However, this expected rigidity is at the same time one of the premises to such a rich theory.

Let $X$ be a smooth complex projective variety. A stable map $f:C \to X$ of class $\beta$ and genus $g$ will be called rigid, if $(C,f)$ is the only point in a connected component of the moduli stack $\overline{\mathcal M}_g(X,\beta)$ of stable maps. The slightly stronger notion of infinitesimal rigidity amounts to $f$ having only trivial infinitesimal deformations. The naive expectation is that such rigid stable maps occur when
$$ \text{vdim } \overline{\mathcal M}_g(X,\beta) = \int_\beta c_1(X) + (g-1)(3-\dim X) = 0, $$
so we expect an abundance of such rigid stable maps if and only if $X$ is a Calabi-Yau threefold. Of course, if $\text{vdim } \overline{\mathcal M}_g(X,\beta) > 0$, then rigid stable maps cannot exist, whereas if $\text{vdim } \overline{\mathcal M}_g(X,\beta) < 0$, then at the very least they must be obstructed. 

Although the questions concerning rigidity can in a sense be circumvented thanks to the excess intersection machinery of Gromov-Witten theory, many such problems, most notably, the Clemens conjecture \cite{[Cl86], [JK96], [Co05], [Co12]}, remain wide open. For generic quintics, Clemens proved the existence of smooth infinitesimally rigid rational curves of arbitrarily large degree as an ingredient in his theorem on algebraic versus homological equivalence of 1-cycles \cite{[Cl83]} and later Katz strengthened the argument \cite{[Ka86a]} to show that such rational curves occur in all degrees. This approach was taken further by Knutsen to prove the existence of smooth infinitesimally rigid curves of genus $g \leq 22$ \cite{[Kn12]}, if $d$ is sufficiently large. 

The main application of the construction presented in the rest of the paper concerns the existence of rigid stable maps of arbitrary genus. Nevertheless, much of the paper is guided by the rational case, which is the best suited to our construction. 

\begin{thm}
Let $g \geq 0$, $d \gg_g 0$ and $Q$ a general quintic threefold. Then $Q$ admits a rigid, degree $d$ stable map $f: C \to Q$ whose source $C$ is a smooth curve of genus $g$.
\end{thm}

The theorem above does not rule out the possibility of $f$ factoring through an unramified cover of smooth curves with source $C$, but this situation can only occur for special pairs $(d,g)$. Refining only the combinatorial part of the arguments below should fix this shortcoming and also prove infinitesimal rigidity, but we have chosen not to pursue these problems here.

\subsection{The approach by specialization} 

The idea that one could gain insight into questions about rational curves (or more generally, arbitrary algebraic curves) on Calabi-Yau threefolds by degenerating the underlying threefold has been around for a long time. Folklore states that the idea originated as a potential approach to the Clemens conjecture. Very naively, if one could exhibit a sufficiently general family of quintics specializing to a reducible variety, locate the flat limits of the curves on the degenerate variety and prove they are finitely many, finiteness would obviously follow. In enumerative geometry, the Gross-Siebert program \cite{[GS06], [GS10]} aims to use toric degenerations and tropical geometry to understand mirror symmetry.

A simple but already very interesting example is the degeneration of a family of quintics to the reducible hypersurface $X_0X_1X_2X_3X_4 = 0$, the union of all coordinate hyperplanes. The project of describing explicitly the limits of the rational curves from the general fibers has been carried out in degree $d=1$, see \cite{[Ka86b], [Ni15]}. However, extending this analysis to arbitrary degree is an extremely difficult task. 

The main purpose of this paper is to exhibit an alternative degeneration of quintic threefolds which is not toric, but provides a different type of leverage in the hunt for the limits of the rational curves. In particular, it is straightforward to give a complete description of all genus zero admissible stable maps to the degenerate fiber, which hasn't been previously achieved with any other degenerations of Calabi-Yau threefolds.

The degeneration we will consider is obtained by running semistable reduction on a sufficiently general pencil of quintic threefolds specializing to the chordal variety of a normal elliptic curve $E$ in ${\mathbb P}^4$. We will obtain a family $W \to \Delta = {\mathbb A}^1$ with $W_0 = Y_1 \cup Y_2$ having only two components. The degeneration has the usual desirable properties (simple normal crossing components in the central fiber, smooth total space etc.) and, most importantly, the two components of the degenerate threefold admit useful maps $\varphi_i: Y_i \to E$. Furthermore, $Y_1 \cap Y_2 \cong E \times E$ and the two maps $\varphi_i$ restrict on this intersection to projections to the two factors. The picture will be described in great detail in section 2. While the analysis itself is uninspiring, the property above should justify this exercise in semi-stable reduction.

\subsection{Limits of algebraic curves} 

It seems useful to illustrate the analysis of limit curves in $W_0$ which will be carried out in section 3 in the case $g=0$. The fact that both components admit maps to genus one curves confines all the irreducible components of any limit to the fibers of these maps. Consider the function
\begin{center}
\begin{tikzpicture}
\node at (-3.5,0) {$\left\{ \begin{matrix} 1\text{-cycles on }W_0 \text{ which are limits} \\ \text{of rational curves from nearby fibers} \\ \end{matrix} \right\}$};
\node at (3.5,0) {$\left\{ \begin{matrix} \text{finite collections of} \\ \text{closed points on }E \times E \\ \end{matrix} \right\}$};
\draw[-{Stealth[length=2mm]}, thin] (0,0) -- (1.1,0);
\node at (0.5,0.3) {$\zeta$};
\end{tikzpicture}
\end{center}
which associates to each limit-$1$-cycle, the points of intersection with the double surface $Y_1 \cap Y_2 \cong E \times E$. Denote by $\zeta_d$ the restriction of $\zeta$ to the subset of limits of degree $d$ curves.

We want to think of the problem of finding the limits in two steps: (1) first, determine the points of intersection with $E \times E$; then (2) for a given such collection of points, find the curves suiting them. It is easy to see that $\zeta_d$ is finite-to-one, which amounts to a satisfactory answer to the second step. Regarding (1), we discover that the patterns of points are related to some intricate combinatorics involving ${\mathbf Q}/{\mathbf Z}$-valued linear algebra on the dual graph. 

The biggest issue is the following: these combinatorial constrains only account for the "kissing condition" and degenerate curves which aren't really limits in the classical sense may occasionally creep in, as in examples 3.5 and 3.6. Therefore, the natural thing to do is to describe the genus zero degenerate stable maps to $W_0$ and its expansions, as defined by Jun Li \cite{[Li01], [Li02]}. Moreover, we will see that nothing essential changes for $g>0$ if we instead impose the condition that all components of the source are rational -- or more accurately, the arithmetic genera of the connected components of the glued relative stable maps are zero. 

The partial extension to $g \geq 0$ is easily sufficient to prove Theorem 1.1. The existence of a perfect obstruction theory ensures that rigid degenerate stable maps smooth out to rigid stable maps and from this point, everything can easily be controlled by purely combinatorial means. After several reductions, we boil down the proof to the existence of graphs with $\smash{\text{rank }\mathrm{H}^1(G,{\mathbf Z})=g}$ and $d$ edges, having the (normalized Laplace) eigenvalue $5/3$ in some ${\mathbf F}_p$ but not in $\smash{\overline{\mathbf Q}}$ and satisfying some additional properties. This final combinatorial step was computer-aided in the following sense: the "building blocks" exhibited in Claim 3.10 were found using a computer, but their validity can also be ascertained by hand with sufficient patience, and their existence is completely natural.  

\bigskip

\noindent \textit{Acknowledgments.} In writing this paper, I have benefited greatly from long and very interesting discussions with Gabriel Bujokas, Xi Chen and my dissertation adviser, Joe Harris. One chapter of my thesis is based on this work. I would also like to thank Herb Clemens, Jun Li, Ravi Vakil and several of my former colleagues at Harvard University for valuable comments and insights. 

Finally, I would like to thank the referee for an extremely careful reading of the submitted manuscript and many valuable suggestions.

\section{Explicit description of the degeneration}

Let $E \subset {\mathbb P}^4$ be a smooth nondegenerate curve of genus $1$ and degree $5$ and let $\Theta$ be the chordal variety of $E$, that is, the singular threefold swept out by all lines connecting pairs of point on $E$. Consider the diagram
\begin{center}
\begin{tikzpicture}

\matrix [column sep  = 10mm, row sep = 10mm] {
	\node (nw) {$\tilde{\Theta}$};
	& \node (ne) {${\mathscr U}$}; 
	& \node (nee) {${\mathbb P}^4$}; \\
	\node(sw) {$E^{(2)}$}; &
	\node (se) {${\mathbb G}(1,4)$}; \\
};
\draw[-{Stealth[length=2mm]}, thin] (ne) -- (nee);
\draw[-{Stealth[length=2mm]}, thin] (nw) -- (sw);
\draw[-{Stealth[length=2mm]}, thin] (nw) -- (ne);
\draw[-{Stealth[length=2mm]}, thin] (sw) -- (se);
\draw[-{Stealth[length=2mm]}, thin] (ne) -- (se);

\end{tikzpicture}
\end{center}
where $E^{(k)}$ is the $k$th symmetric power of $E$ and ${\mathscr U} \rightarrow {\mathbb G}(1,4)$ is the universal family of lines in ${\mathbb P}^4$. The map $E^{(2)} \to {\mathbb G}(1,4)$ is defined by sending a length two subscheme of $E$ to the unique line in ${\mathbb P}^4$ which contains it. Define $\smash{ \tilde{\Theta} = E^{(2)} \times_{{\mathbb G}(1,4)} {\mathscr U} }$ and $\Theta$ its image in ${\mathbb P}^4$. It is not hard to check that $\Theta$ has degree $5$. We'll briefly review this calculation and other easy geometric properties early in \S2.1. 

Consider a general pencil of quintic threefolds whose fiber over $0 \in {\mathbb A}^1 \subset {\mathbb P}^1$ is the secant variety $\Theta$. The total space of the pencil is the variety $\{F_\Theta + t F=0\} \subset {\mathbb A}_t^1 \times {\mathbb P}^4$ where $F_\Theta$ is the homogeneous polynomial of $\Theta$ and $F \in {\mathrm{H}}^0({\mathscr O}(5))$. Seen over ${\mathbb A}^1_t$, this is a family of quintics degenerating to $\Theta$. The main observation is that running semistable reduction on this family yields a central fiber with two (smooth and transverse) components both admitting nonconstant maps to a smooth algebraic curve of genus one.

The semistable reduction calculation will be carried out in \S2.2. There are two steps: (1) make a base change of order $3$ totally ramified at $0 \in {\mathbb A}^1$; and (2) blow up the curve $\{ 0 \} \times E$ in the total space of the family. Said differently, $W$ is the proper transform of $Z = \{F_\Theta + t^3 F = 0\}$ in the blowup of ${\mathbb A}^1 \times {\mathbb P}^4$ along $\{0\} \times E$. We will check that: (1) the total space $W$ is smooth; (2) its central fiber $W_0$ has two normal crossing smooth components: $\smash{Y_1 = \tilde{\Theta}}$, which is fibered over $\text{Pic}^2 (E)$ with fibers ${\mathbb F}^1$ and the exceptional divisor $Y_2$, which admits a map to $E$, whose fibers are generically smooth "cyclic" cubic surfaces, with the exception of $25$ projective cones over a smooth plane cubic, above the points $p_1,...,p_{25} \in E$ where the base locus of the pencil meets $E$; (3) the intersection $Y_1 \cap Y_2$ is isomorphic to $E \times E$; and (4) the two maps from $Y_i$ to $\mathrm{Pic}^2(E)$ or $E$ restrict on the intersection to the tensor product map $E \times E \to \text{Pic}^2(E)$ and projection to one of the factors respectively.

It may be useful to abstract away the essential features of the configuration above. What we have is a family $\pi:W \to \Delta$ over a smooth curve $\Delta$, with smooth total space, such that $W_0 = Y_1 \cup Y_2$ with $Y_1 \cap Y_2 = E_1 \times E_2$ and $E_1 \cong E_2$ isomorphic smooth genus $1$ curves, $Y_i$ smooth and admitting a regular map $\varphi_i:Y_i \to E_i$ such that $\varphi_i$ restricts on $E_1 \cap E_2$ to projection to the $i$th factor, for $i=1,2$. Moreover, the fibers of the pair $(Y_i,E_1 \times E_2)$ over $E_i$ are generically pairs consisting of a smooth surface and a smooth effective anticanonical divisor. In this paper we will only be concerned with the explicit example exhibited here, but it seems very likely that degenerations fitting into the pattern above also exist for other deformation classes of Calabi-Yau threefolds.

The purely aesthetic change of perspective from $E \times E$ to $E_1 \times E_2$ is obtained by considering the isomorphism $(\otimes, \text{proj}_1): E \times E \to E_1 \times E_2$, where $E_1 = \text{Pic}^2 (E)$ and $E_2 = E$. A choice of one of the $25$ hyperflexes of $E \subset {\mathbb P}^4$ gives a convenient identification of $E_1 = \text{Pic}^2(E)$ with $E$. Of course, if $h$ is the hyperflex, the isomorphism is the map $p \mapsto {\mathscr O}_E(p+h)$.

\subsection{Chordal variety of a normal elliptic curve} 
We begin with an elementary review of the chordal variety of a normal elliptic curve in ${\mathbb P}^4$. We use the same notation as above: $E$ is the normal elliptic curve and $\Theta$ is the chordal variety, which by definition is the singular threefold swept out by all the lines intersecting $E$ in a scheme of length $2$. Let $\tilde{P}$ be the blow up of ${\mathbb P}^4$ along the curve $E$. We have an obvious tower of ${\mathbb P}^1$-bundles  $\tilde{\Theta} \to E^{(2)} \to \text{Pic}^2 (E)$ and thus $\tilde{\Theta}$ is smooth. It follows that for any $D \in \text{Pic}^2 (E)$, the fiber $\tilde{\Theta}_{D}$ is some rational ruled surface ${\mathbb F}_n$. 

\begin{lemma}
The degree of $\Theta$ in ${\mathbb P}^4$ is equal to $5$.
\end{lemma}

\begin{proof}
We may project away from a general line $\ell \subset {\mathbb P}^4$ onto ${\mathbb P}^2$ and note that intersection points of $\ell$ with $\Theta$ correspond to the nodes of the projection of $E$. Since the projection of $E$ to ${\mathbb P}^2$ has degree $5$ and geometric genus $1$, the genus-degree formula implies that the number of nodes is equal to $5$. \end{proof}

\begin{lemma}
The map $j:\smash{\tilde{\Theta} \to {\tilde{P}}}$ is a closed embedding. 
\end{lemma}

\begin{proof} 
It suffices to show that $j$ is injective and $dj:\smash{ T \tilde{\Theta} \to j^* T{\tilde{P}}}$ is fiberwise injective. For the first claim, assume that there exist two distinct points $p,q \in \smash{\tilde{\Theta}}$ such that $j(p) = j(q)$. Let $\smash{\tilde{\ell}_p, \tilde{\ell}_q} \subset \tilde{\Theta}$ be the two secant lines containing $p$ and $q$ and let $\ell_p,\ell_q$ be their images in ${\mathbb P}^4$. Clearly, $\ell_p \neq \ell_q$. Assume first that $\ell_p \cap \ell_q \notin E$. Then the length $4$ subscheme $D_{p,q} = E \cap (\ell_p \cup \ell_q)$ of $E$ is contained in the projective subspace of ${\mathbb P}^4$ spanned by $\ell_p$ and $\ell_q$, which has dimension $2$. However, by Riemann-Roch, 
$$ h^0({\mathscr O}_E(1) \otimes {\mathscr O}_E(-D_{p,q})) = 1, $$
leading to contradiction. If $\ell_p \cap \ell_q \in E$, then $\smash{\tilde{\ell}_p \cap \tilde{\ell}_q = \emptyset}$, which is impossible as well. 

The second claim is an infinitesimal version of the first. Assume that there is a tangent vector of to $\smash{\tilde{\Theta}}$ which is killed by $dj$. Then there exists a line and a first order thickening $\smash{\tilde{\ell} \subset \hat{\ell} \subset {\tilde{\Theta}}}$ such that the corresponding section of $\smash{{\mathscr N}_{\tilde{\ell}/\tilde{P}}}$ vanishes at some point, or equivalently, ${\ell}'$ is contained in a plane, where $\smash{\ell \subset \ell'}$ is the image in ${\mathbb P}^4$ of $\smash{\tilde{\ell} \subset \hat{\ell}}$. Then an argument similar to the one above leads to contradiction.
\end{proof}

\begin{lemma}
All the fibers of $\tilde{\Theta} \to \mathrm{Pic}^2(E)$ are isomorphic to ${\mathbb F}_1$.
\end{lemma}

\begin{proof}
One way to prove this is via an indirect argument starting with the Hirzebruch surface. Let $f, e_{\infty}, e_0 \in \text{Pic} ({\mathbb F}_1)$ be the classes of the fibration, of the directrix and of the preimage of a general line from the projective plane obtained by blowing down the directrix respectively. These classes generate $\text{Pic} ({\mathbb F}_1)$, subject to the relation $e_0 = e_\infty + f$. The line bundle ${\mathscr O}_{{\mathbb F}_1} (e_\infty + 2f)$ defines the standard embedding of ${\mathbb F}_1$ into ${\mathbb P}^4$. 

Consider the anticanonical class  $- K_{{\mathbb F}_1} = {\mathscr O}_{{\mathbb F}_1} (2e_\infty + 3f)$. The divisors in the associated linear system are curves of arithmetic genus one, intersecting the directrix once and the fibers of the ruling twice. It is not hard to see that we can find such a divisor, which is abstractly isomorphic to our curve $E$ and such that the divisor class on $E$ cut by the ruling is precisely $D$. Then, embedding the pair $E \subset {\mathbb F}_1$ in ${\mathbb P}^4$, it is not hard to see that $E$ is an elliptic normal curve and that the surface ${\mathbb F}_1$ is precisely the surface swept out by the lines connecting pairs of points on $E$ whose sum as divisor classes is precisely $D$. Since the action of $\text{PGL}(5, {\mathbb C})$ on the set of normal elliptic curves of a fixed $j$-invariant is transitive, this argument proves the claim that $\tilde{\Theta}_{D}$ is an ${\mathbb F}_1$ surface. 

Note that the lemma also follows from the observation that $\tilde{\Theta}_{D}$ is a minimal degree variety, due to the classification proved in \cite{[EH87]}.
\end{proof}

Note that the intersection $S$ of $\tilde{\Theta}$ with the exceptional divisor $T$ of $\tilde{P}$ is canonically isomorphic to $E \times E$. Roughly, any point in the intersection is specified by choosing a point on $E$ and an infinitesimal direction inside $\Theta$ from the chosen point, the latter also being specified by a point on $E$. The restriction of $\tilde{\Theta} \to \mathrm{Pic}^2 (E)$ to this locus is the tensor product map 
$$ \nu: E \times E \longrightarrow \mathrm{Pic}^2 (E) $$
by construction. In this section, we order the factors of $S = E \times E$ such that the restriction $S \to E$ of the blowdown map is projection to the first factor. We will usually identify $\nu^{-1}(D)$ with $E$ in this way. In the statement below, $S_p$ and $T_p$ denote the fibers of the maps $S \to E$ and $T \to E$ over the point $p \in E$.

\begin{lemma} 
$(a)$ The point where the directrix of $\tilde{\Theta}_{D}$ intersects $\nu^{-1}({D}) \cong E$ is the one whose divisor class is ${\mathscr O}_E(1) \otimes {\mathscr O}_E(-2D)$.

$(b)$ The restriction of the hyperplane class of $T_p$ to $S_p \cong E$ is ${\mathscr O}_E(1) \otimes {\mathscr O}_E(-2p)$, for all closed $p \in E$.
\end{lemma}

\begin{proof}
Part $(a)$ follows from the proof of Lemma 2.3. For part $(b)$, simply pick a hyperplane containing the line tangent to $E \subset {\mathbb P}^4$ at $p$ and project away from the tangent line. Note that $S_p$ sits inside $T_p \cong {\mathbb P}^2$ as a smooth plane cubic.
\end{proof}

Denote the blowdown map $\tilde{P} \to {\mathbb P}^4$ by $\rho$ and its restriction to $T$ by 
$$ \tau = \rho|_T: T \longrightarrow E. $$
Of course, $T$ is the projectivization of the normal bundle of $E$ and as such comes equipped with a tautological line bundle ${\mathscr O}_T(-1)$ isomorphic to ${\mathscr O}_{\tilde{P}}(T) \otimes {\mathscr O}_T$.

\begin{claim}
We have $\smash{ \rho^* {\mathscr O}_{{\mathbb P}^4}(\Theta)  \cong {\mathscr O}_{\tilde{P}}(\tilde{\Theta} + 3T) }$ and $\smash{ {\mathscr O}_T(S) \cong {\mathscr O}_T(3) \otimes \tau^* {\mathscr O}_E(5) }$.
\end{claim}

\begin{proof}
We must have $\rho^* \Theta = \tilde{\Theta} + k T$ for some $k$ as Cartier divisors. For a line $L$ contained inside the fibers of $\tau$, we have $(\rho^{*} \Theta \cdot L) = 0$,  $(\tilde{\Theta} \cdot L)=3$ since $S_p \subset T_p$ has degree $3$ and $(T \cdot L) = -1$, hence $k=3$. For the second part, restrict to $T$. The left hand side and the right hand side become respectively $\smash{ \rho^*{\mathscr O}_{{\mathbb P}^4}(5) \otimes {\mathscr O}_T \cong \tau^* {\mathscr O}_E(5) }$ and $\smash{ {\mathscr O}_{\tilde{P}}(\tilde{\Theta} + 3T) \otimes {\mathscr O}_T \cong {\mathscr O}_T(-3) \otimes {\mathscr O}_T(S) }$, proving the claim. 
\end{proof}

\begin{lemma}
The ideal sheaf ${\mathscr I}_{\Theta/{\mathbb P}^4}$, seen inside ${\mathscr O}_{{\mathbb P}^4}$, is contained in ${\mathscr I}_{E/{\mathbb P}^4}^3$.
\end{lemma}

\begin{proof} 
From the claim, we have a map of ${\mathscr O}_{\tilde P}$-modules $\smash{ {\mathscr O}_{\tilde{P}}(3T) \to \rho^* {\mathscr O}_{{\mathbb P}^4}(\Theta) }$. Let $\smash { \rho^*{\mathscr I}_{{\Theta/{\mathbb P}^4}} \to {\mathscr I}^3_{T/\tilde{P}} }$ be the dual. We obtain a map $\smash{ {\mathscr I}_{{\Theta/{\mathbb P}^4}} \to \rho_*{\mathscr I}^3_{T/\tilde{P}} }$ from the adjoint property and it is not hard to check that the latter sheaf is simply $\smash{ {\mathscr I}_{E/{\mathbb P}^4}^3 }$.
\end{proof}

\begin{rem}
Returning to the discussion in the introduction concerning the degeneration to a union of five hyperplanes, we should point out that a union of five hyperplanes is in fact a limit of a family of such chordal varieties. Indeed, let $\ell_{i,j}$ be the line $\{ [\alpha e_i + \beta e_j] : \alpha,\beta \in {\mathbf C} \}$. Then $\ell_{0,1} \cup \ell_{1,2} \cup \ell_{2,3} \cup \ell_{3,4} \cup \ell_{4,0}$ is a reducible curve of degree $5$ and arithmetic genus $1$ whose chordal variety is simply the union of the coordinate hyperplanes. It is not hard to check that the curve is smoothable.
\end{rem}
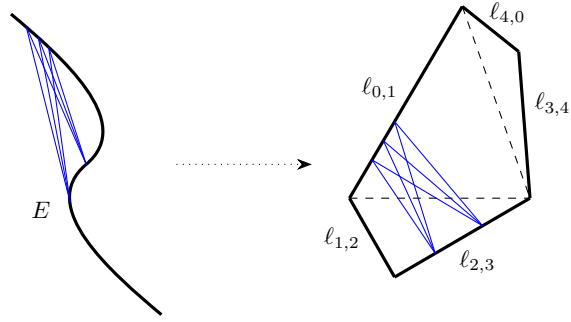
\begin{center}
\begin{tikzpicture}[every node/.style={scale=0.9}]

\def\z{2.5}
\def\k{1.5}
\def\ax{0}
\def\ay{1.4*\k}
\def\bx{-1*\k}
\def\by{-0.3*\k}
\def\cx{-0.6*\k}
\def\cy{-1*\k}
\def\dx{0.6*\k}
\def\dy{-0.3*\k}
\def\ex{0.5*\k}
\def\ey{1*\k}

\def\ku{0.2}
\def\kv{0.3}
\def\kw{0.4}
\def\kt{0.35}
\def\ks{0.7}

%\node at (\z,0.3) {$Z_4 = 0$};
\node at (\z-1.1,1) {$\ell_{0,1}$};
\node at (\z-1.6,-1) {$\ell_{1,2}$};
\node at (\z+0.2,-1.3) {$\ell_{2,3}$};
\node at (\z+1.2,0.8) {$\ell_{3,4}$};
\node at (\z+0.6,2) {$\ell_{4,0}$};
\node at (-\z-0.6,-0.6) {$E$};

%\draw(\ku,2) -- (\ks,2); 

\def\ux{\ku * \ax + \bx -\ku * \bx + \z}
\def\uy{\ku * \ay + \by -\ku * \by}

%\draw (1,1) -- (\ux,\uy);

\def\vx{\kv * \ax + \bx - \kv * \bx + \z}
\def\vy{\kv * \ay + \by - \kv * \by}

\def\wx{\kw * \ax + \bx - \kw * \bx + \z}
\def\wy{\kw * \ay + \by - \kw * \by}

\def\tx{\kt * \cx + \dx - \kt * \dx + \z}
\def\ty{\kt * \cy + \dy - \kt * \dy}

\def\sx{\ks * \cx + \dx -\ks * \dx + \z}
\def\sy{\ks * \cy + \dy - \ks * \dy}

\draw[-{Stealth[length=2mm]}, dotted] (-1.3,0) -- (0.5,0);

\draw [very thick] (-\z,0) to [in= 140, out = -140]++(1,-2);
\draw [very thick] (-\z,0) to [in= -40, out = +40]++(-1,2);

\draw [very thick] (\ax+\z,\ay) -- (\bx+\z,\by);
%\draw [dashed] (\ax+\z,\ay) -- (\cx+\z,\cy);
\draw [dashed] (\ax+\z,\ay) -- (\dx+\z,\dy);
\draw [very thick] (\ax+\z,\ay) -- (\ex+\z,\ey);
\draw [very thick] (\bx+\z,\by) -- (\cx+\z,\cy);
\draw [dashed] (\bx+\z,\by) -- (\dx+\z,\dy);
%\draw [dashed] (\bx+\z,\by) -- (\ex+\z,\ey);
\draw [very thick] (\cx+\z,\cy) -- (\dx+\z,\dy);
%\draw [dashed] (\cx+\z,\cy) -- (\ex+\z,\ey);
\draw [very thick] (\dx+\z,\dy) -- (\ex+\z,\ey);

\draw [blue] (\ux,\uy) -- (\tx,\ty);
\draw [blue](\vx,\vy) -- (\tx,\ty);
\draw [blue](\wx,\wy) -- (\tx,\ty);
\draw [blue](\ux,\uy) -- (\sx,\sy);
\draw [blue](\vx,\vy) -- (\sx,\sy);
\draw [blue](\wx,\wy) -- (\sx,\sy);

\draw [blue](-0.215-\z,-0.5) -- (-0.5-\z,1.56);
\draw [blue](-0.215-\z,-0.5) -- (-0.64-\z,1.7);
\draw [blue](-0.215-\z,-0.5) -- (-0.8-\z,1.84);
\draw [blue](-\z,0) -- (-0.5-\z,1.56);
\draw [blue](-\z,0) -- (-0.64-\z,1.7);
\draw [blue](-\z,0) -- (-0.8-\z,1.84);
\end{tikzpicture}
\end{center}

\subsection{Base change and blowup} 
The previous lemma suggests that in order to obtain a good central fiber, it is necessary to perform a base change of order $3$ before blowing up the singular locus. Set ${\mathscr N} = {\mathscr N}_{E/{\mathbb P}^4}$. The normal bundle of $E \times \{0\}$ in ${\mathbb P}^4 \times {\mathbb A}^1$ is isomorphic to ${\mathscr N} \oplus {\mathscr O}_E$. To avoid being confusing later, we should clarify that the summand ${\mathscr O}_E$ should actually be written as ${\mathscr O}_E \langle t^{-1} \rangle$, where $t$ is the affine coordinate of ${\mathbb A^1}$. 

Let $X$ be the blowup of ${\mathbb P}^4 \times {\mathbb A}^1$ along $E \times \{0\}$. From now on we will usually write just $E$ instead of $E \times \{0\}$. The central fiber $W_0$ has two components: $X_{0,1} \cong \tilde{P}$ and $ X_{0,2} = \mathrm{\mathrm{Proj} Sym } ({\mathscr N} \oplus {\mathscr O}_E)^{\vee}$.
The inclusion ${\mathscr O}_E \subset {\mathscr O}_E \oplus {\mathscr N}$ induces a section 
$$ \sigma:E \longrightarrow X_{0,2}, $$ 
whose image is geometrically the intersection of the proper transform of $E \times {\mathbb A}^1$ with the central fiber. Dually, it defines a distinguished section $t$ of ${\mathscr O}_{X_{0,2}}(T)$ which we may identify with the affine coordinate above. 

After a base change of order $3$, a pencil of quintics specializing to $\Theta$ has total space $Z = \{F_\Theta + t^3 F=0\}$ where we choose $F$ to be a very general homogeneous quintic polynomial and $F_\Theta$ is the homogeneous equation of $\Theta$. The central fiber $W_0$ of the proper transform $W$ of $Z$ correspondingly has two components: $Y_1 = \tilde{\Theta} \subset X_{0,1}$ and $Y_2 \subset X_{0,2}$. There seems to be little choice in describing $Y_2$ besides writing explicit equations. Of course, by definition, $Y_2 = \text{\text{Proj}}_E \oplus_{k \geq 0} {\mathscr I}_{E/Z}^k/{\mathscr I}_{E/Z}^{k+1} $. We will show that $Y_2$ is the total space of a family of cubic surfaces in the fibers of the projective bundle 
$$ \varphi_2:X_{0,2} \longrightarrow E. $$
Note that Lemma 2.6 further implies that 
$$ {\mathscr I}_{Z/{\mathbb P}^4 \times {\mathbb A}^1} \subset {\mathscr I}^3_{E/{\mathbb P}^4 \times {\mathbb A}^1}. \eqno(2.1) $$
Indeed, given the shape of the equation defining $Z$, we have
$$ {\mathscr I}_{Z/{\mathbb P}^4 \times {\mathbb A}^1} \subseteq {\mathscr I}_{\Theta/{\mathbb P}^4} \boxtimes {\mathscr O}_{{\mathbb A}^1} + {\mathscr I}^3_{{\mathbb P}^4 \times \{0\}/{\mathbb P}^4 \times {\mathbb A}^1} \subseteq  {\mathscr I}^3_{E/{\mathbb P}^4 \times {\mathbb A}^1}. $$
Consider the following commutative diagram of coherent sheaves on $E$:
\begin{center}
\begin{tikzpicture}

\matrix [column sep  = 12mm, row sep = 10mm] {
	\node (nw) {${\mathscr I}_{Z/{\mathbb P}^4 \times {\mathbb A}^1} \otimes {\mathscr O}_E$}; &
	\node (nc) {${\mathscr I}^3_{E/{\mathbb P}^4 \times {\mathbb A}^1} \otimes {\mathscr O}_E$}; &
	\node (ne) {$\text{Sym}^3 ({\mathscr N} \oplus {\mathscr O}_E)^\vee$}; \\
	\node (sw) {${\mathscr I}_{\Theta/{\mathbb P}^4} \otimes {\mathscr O}_E$}; &
	\node (sc) {${\mathscr I}^3_{E/{\mathbb P}^4} \otimes {\mathscr O}_E$}; &
	\node (se) {$\text{Sym}^3 {\mathscr N}^\vee$}; \\
};

\draw[-{Stealth[length=2mm]}, thin] (nw) -- (nc);
\draw[-{Stealth[length=2mm]}, thin] (nc) -- (ne);
\draw[-{Stealth[length=2mm]}, thin] (sw) -- (sc);
\draw[-{Stealth[length=2mm]}, thin] (sc) -- (se);
\draw[-{Stealth[length=2mm]}, thin] (nw) -- (sw);
\draw[-{Stealth[length=2mm]}, thin] (nc) -- (sc);
\draw[-{Stealth[length=2mm]}, thin] (ne) -- (se);

\end{tikzpicture}
\end{center}
The left vertical arrow is an isomorphism as it is simply the restriction to $E$ of the isomorphism ${\mathscr I}_{Z/{\mathbb P}^4 \times {\mathbb A}^1}|_{\{t=0\}} \cong {\mathscr I}_{\Theta/{\mathbb P}^4}$. The horizontal arrows on the right side are also isomorphisms by well known properties. Taking into account all these isomorphisms, the diagram above becomes the following diagram:
\begin{center}
\begin{tikzpicture}

\matrix [column sep  = 12mm, row sep = 5mm] {
	& \node (ne) {$\varphi_{2*}{\mathscr O}_{X_{0,2}}(3)$}; & \\
	\node (w) {${\mathscr O}_E(-5)$}; & \\
	& \node (se) {$\tau_{*}{\mathscr O}_{T}(3)$}; & \\	
};

\draw[-{Stealth[length=2mm]}, thin] (w) -- (ne);
\draw[-{Stealth[length=2mm]}, thin] (w) -- (se);
\draw[-{Stealth[length=2mm]}, thin] (ne) -- (se);

\end{tikzpicture}
\end{center}
By push-pull, we get a global section $y_2$ of ${\mathscr O}_{X_{0,2}}(3) \otimes \varphi_2^*{\mathscr O}_E(5)$,  which restricts on $T$ to the global section $s \in \text{H}^0({\mathscr O}_{T}(3) \otimes \tau^*{\mathscr O}_E(5))$ defining $S$ inside $T$, due to the second part of Claim 2.5.

\begin{claim}
The scheme-theoretic vanishing locus of $y_2$ is precisely $Y_2$. Furthermore, there exists a natural identification
$$ \text{H}^0({\mathscr O}_{X_{0,2}}(3) \otimes \varphi_2^*{\mathscr O}_E(5)) \cong \bigoplus_{k=0}^{3} \text{H}^0({\mathscr O}_{T}(k) \otimes \tau^*{\mathscr O}_E(5)) \otimes {\mathbb C}\langle t^{3-k}\rangle, \eqno(2.2) $$
under which
$$ y_2 = \tau^* F|_E \otimes t^3 +  s \otimes 1, \eqno(2.3) $$
where $F|_E$ is the restriction of $F$ to $E$ and $\tau^*F|_E$ is its pullback to $T$.
\end{claim}

\begin{proof}
The fact that $y_2$ vanishes on $Y_2$ can be seen from construction as follows. By (2.1), $\smash{{\mathscr I}_{Z/{\mathbb P}^4 \times {\mathbb A}^1} \otimes {\mathscr O}_E}$ lies in the kernel of the map $\smash{{\mathscr I}^3_{E/{\mathbb P}^4 \times {\mathbb A}^1}|_E \to {\mathscr I}^3_{E/Z}|_E}$, which is a piece of the map of graded ${\mathscr O}_E$-algebras
$$ \bigoplus_{k \geq 0} {\mathscr I}^k_{E/{\mathbb P}^4 \times {\mathbb A}^1} \otimes {\mathscr O}_E \longrightarrow \bigoplus_{k \geq 0} {\mathscr I}^k_{E/Z} \otimes {\mathscr O}_E $$
corresponding to the closed immersion $Y_2 \to X_{0,2}$, thus proving that $y_2$ vanishes on $Y_2$. The equation defining $Y_2$ scheme-theoretically is $y_2$, since $y_2$ restricts on $T$ to the section $s$ defining $S$.

For the purpose of giving explicit equations, it is most convenient to use the description of our sheaves as symmetric powers of normal bundles. Standard multilinear algebra gives a canonical decomposition
$$ \text{Sym}^3 ({\mathscr N} \oplus {\mathscr O}_E \langle t^{-1} \rangle)^\vee \cong \bigoplus_{k=0}^{3} \text{Sym}^k {\mathscr N}^\vee \otimes \langle t^{3-k}\rangle {\mathscr O}_E, $$
or equivalently $\varphi_{2*}{\mathscr O}_{X_{0,2}}(3) \cong \oplus_{k=0}^3 \tau_{*}{\mathscr O}_{T}(k)\otimes \langle t^{3-k}\rangle {\mathscr O}_E$, which can be rearranged using push-pull once more as (2.2) above. The explicit formula (2.3) for $y_2$ follows simply by unwinding all the definitions.
\end{proof}

Let $p_1,p_2,...,p_{25}$ be the points of intersection of $E$ with the base locus of the pencil of quintic threefolds.

\begin{proposition}
The scheme $Y_2$ is smooth and irreducible. For closed $p \in E$, the fiber of $\varphi_2:Y_2 \to E$ over $p$ is a cubic surface inside the corresponding fiber of $W_{0,2}$ which is either a cone over $E$ with vertex at $\sigma(p)$ if $p \in \{p_1,p_2,...,p_{25}\}$ or a special smooth cubic surface which can be expressed as a triple cover of $T_p$ totally ramified at $S_p \cong E \subset T_p$ via projection from $\sigma(p)$ otherwise.
\end{proposition}

\begin{proof}
Most statements follow immediately from $(2.3)$ since the formula shows that the equation of each fiber of $Y_2$ inside the corresponding fiber of $X_{0,2}$ is of the form $c_pT^3 + F_3(X,Y,Z) = 0$, for some constant $c_p$ which is zero precisely when $p \in \{p_1,p_2,...,p_{25}\}$. The only somewhat nontrivial thing that's left to check is the smoothness of $Y_2$ at the $25$ points $\sigma(p_i)$, but this follows easily from the fact that $\partial c / \partial p(p_i) \neq 0$, which in turn follows from the fact that $p_1,p_2,...,p_{25}$ are distinct. The derivative above is of course well-defined since sections of line bundles admit derivatives at vanishing points.
\end{proof}
 
\begin{proposition} 
For a general choice of the pencil of quintics containing $\Theta$, $Y_1$ and $Y_2$ are smooth irreducible divisors on $W$ meeting transversally and $W$ is smooth in a neighborhood of the central fiber.
\end{proposition}

\begin{proof}
It's worth noting that, since $W$ is a Cartier divisor on $X$, there will be no need to worry about embedded components in this proof. Smoothness and irreducibility of $Y_1$ and $Y_2$ have already been checked. 

First, let's verify smoothness of the total space. It suffices to check smoothness at closed $p \in W_0$, say $p \in Y_i$. The idea now is to use once more the embedding in $X$. If $W$ were singular at $p$, then the intersection $Y_i = W \cap X_{0,i}$ would also be singular at $p$, which is a contradiction. Slightly more formally, $T_pY_i = T_pW \times_{T_pX} T_pX_{0,i}$ and since $T_pY_i \neq T_pX_{0,i}$ then also $T_pW \neq T_pX$ proving smoothness. Transversality of $Y_1$ and $Y_2$ is automatic, since we already know that their intersection is $S$, which is smooth.
\end{proof}

\subsection{Triple cyclic covers of the projective plane} 
We conclude this section with some elementary remarks on the special class of cubic surfaces encountered earlier. If $S$ is a cubic surface whose equation in ${\mathbb P}^3$ is given in coordinates by
$$ F_3(X,Y,Z) + T^3 = 0 $$
for some homogeneous degree $3$ polynomial $F_3 \in {\mathbf C}[X,Y,Z]$ which defines a smooth plane cubic $E$, then projection from the point $[0:0:0:1]$ to the plane $\{T=0\}$ exhibits $S$ as a triple cover totally branched along $E$ and \'{e}tale elsewhere.

As it is the case with all smooth cubic surfaces, $\text{Pic}(S)$ is (over-)generated by the divisor classes of the $27$ lines. However, the configuration of lines is special: the $27$ lines come in $9$ triples of coplanar lines passing through each of the $9$ flex points of $E$. Indeed, the plane determined by the center of projection and any of the triple tangents to $E$ cuts $S$ along a plane cubic with a triple point at the corresponding flex of $E$ and it therefore has to be a union of $3$ lines. Consequently, the image of the restriction map $\text{Pic}(S) \to \text{Pic}(E)$ consists precisely of the line bundles on $E$ whose cube is some ${\mathscr O}_E(k)$. Finally, we prove some elementary results which will be used in the proof of the existence of rigid stable maps of arbitrary genus.

\begin{lemma}
Let $f_1$ and $f_2$ be two distinct flex points of $E$ and $\ell_1$ an arbitrary line in $S$ passing through $f_1$. Then there exists a unique line $\ell_2 \subset S$ passing through $f_2$ which intersects $\ell_1$.  
\end{lemma}

\begin{proof}
We will give an indirect argument proving uniqueness first. If there were two such lines $\ell_2$ and $\ell'_2$, then $\ell_1$ would have to lie in the plane spanned by $\ell_2$ and $\ell'_2$. However, the only other line on $S$ contained in the plane spanned by $\ell_2$ and $\ell'_2$ is just the third line on $S$ passing through $f_2$, which is different from $\ell_1$ thus proving uniqueness. 

Existence can be inferred from the fact that $\ell_1$ has to intersect $10$ other lines on $S$. Two of them are the other lines passing through $f_1$, so there are $8$ more lines to account for. However, there are precisely $8$ flex points of $E$ besides $f_1$, so there can only be precisely one line intersecting $\ell_1$ through each such point.
\end{proof}

\begin{lemma}
Let $k \leq 3$ and $D_E$ be an effective divisor on $E$ such that ${\mathscr O}_E(3D_E) \cong {\mathscr O}_E(k)$. Then there exists an effective divisor $D_S$ on $S$ such that $D_S \cap E = D_E$ and the complete linear system $|D_S|$ is either: (1) a singleton consisting of a line if $k=1$; (2) a pencil of conics if $k=2$; (3) either (3a) a net of twisted cubics or (3b) the anticanonical linear system, if $k=3$. Moreover, in (3b), $D_S$ can be chosen to be singular.
\end{lemma}

\begin{proof}
If $k=1$, we may just pick any of the $3$ lines through $D_E=p$. If $k=2$, there exist two distinct flex points of $E$ $p$ and $q$ such that $D_E \sim p+q$. Let $\ell_p$ be any line through $p$. By Lemma 2.11, there exists a line $\ell_q$ through $q$ which intersects $\ell_p$. Then $|\ell_p+\ell_q|$ is a pencil of conics and it is easy to see that the restriction $|\ell_p+\ell_q| \to |D_E|$ is bijective, so there exists $D_S \in |\ell_p + \ell_q|$ such that $D_S \cap E = D_E$. 

Finally, if $k=3$, there exist $3$ distinct flex points of $E$ $p,q,r$ such that $D_E \sim p+q+r$. Let $\ell_p$ be a line through $p$ and $\ell_q,\ell_r$ the lines passing through $q$ and $r$ which intersect $\ell_p$. Then there are two cases:

(3a) $\ell_q \cap \ell_r = \emptyset$. In this case, $|\ell_p+\ell_q+\ell_r|$ is a net of twisted cubics and the restriction $|\ell_p+\ell_q+\ell_r| \to |D_E|$ is again bijective, so there exists $D_S \in |\ell_p + \ell_q + \ell_r|$ such that $D_S \cap E = D_E$.

(3b) $\ell_q \cap \ell_r \neq \emptyset$. Now $| \ell_p + \ell_q + \ell_r|$ is the anticanonical linear system, i.e. the linear system of hyperplane sections. In this case, the (now rational) restriction map $|\ell_p + \ell_q + \ell_r| \to |D_E|$ is surjective with fibers of dimension one. The fiber of $D_E$ is a pencil of hyperplane sections, so $D_S$ can be any singular member of this pencil.
\end{proof}

\section{Degenerate stable maps to the central fiber}

As emphasized in the introduction, the motivation for the choice of the specific degeneration worked out in the previous section is the presence of the two nonconstant maps to $E$, which greatly restricts the position of the limits of the nearby rational curves and, to a lower extent, of nearby curves of arbitrary genus. The theoretical framework for working with stable maps to degenerations has been laid out in \cite{[Li01], [Li02]}. The degeneration formula for Gromov-Witten invariants proved in loc. cit. will be less important to us than the construction of the moduli spaces, namely the extension of the family of moduli spaces of stable maps over the singular fiber. Very roughly, this degenerate moduli space parametrizes locally deformable (to nearby fibers) maps to $W_0$ with semistable source, with a twist -- namely, the idea to allow an expansion of $W_0$ by inserting a chain of ruled varieties between $Y_1$ and $Y_2$ with the purpose of avoiding contracted components (of the source curve) to the singular locus of $W_0$.  

Recall that to the total space $W$ and the pair $Y_i^{\text{rel}} = (Y_i,S)$ we may associate the Artin stack of expanded degenerations ${\mathfrak W}$, respectively the stack of expanded pairs ${\mathfrak Y}_i^{\text{rel}}$ \cite{[Li01]}. The stacks of expanded degenerations and expanded pairs allow one to define the spaces of stable maps to expanded degenerations and relative stable maps, which will be denoted by ${\mathcal M}({\mathfrak W},\Gamma)$ and ${\mathcal M}({\mathfrak Y}_i^\text{rel},\Gamma_i)$ respectively. The topological data $\Gamma$ consists simply of the triple $(g,d,k=0)$, where $g$ is the arithmetic genus of the source, $d$ is the degree of the stable maps relative to some predetermined relatively ample line bundle ${\mathscr H}$ on $W$ and $k$ is the number of ordinary marked points. The topological data specified in $\Gamma_i$ is the following: 

(T1) a decorated graph with vertices $V(\Gamma_i)$ corresponding to the connected components of the stable maps, roots $R(\Gamma_i)$ corresponding to the $r$ distinguished marked points attached to the suitable vertices and no edges;

(T2) a function $\mu: R(\Gamma_i) \to {\mathbf Z}^+$ assigning the suitable order of intersection with $S$ at each distinguished marked point;

(T3) functions $\beta_i: V(\Gamma_i) \to \text{H}_2(Y_i,{\mathbf Z})$ and $g_i:V(\Gamma_i) \to {\mathbf N}$ prescribing the homology class in $Y_i$, respectively the arithmetic genus of each connected component.

The total space forms a family $\pi^{\mathcal M} : {\mathcal M}({\mathfrak W},\Gamma) \to \Delta = {\mathbb A}^1$ which is proper over $\Delta$ by \cite{[Li01]}. There exist distinguished evaluation morphisms 
$$ {\mathbf q}_i:{\mathcal M}({\mathfrak Y}_i^\text{rel},\Gamma_i) \longrightarrow S^r, $$  
where the distinguished marked points get artificially ordered as part of the topological data $\Gamma_i$, for $i=1,2$. Whenever $\Gamma_1$ and $\Gamma_2$ are compatible in the sense of having the same number $r$ of roots and the corresponding roots are weighted identically, their fiber product admits a morphism
$$ \Phi_{\Gamma}: {\mathcal M}({\mathfrak Y}_1^\text{rel},\Gamma_1) \times_{S^r} {\mathcal M}({\mathfrak Y}_2^\text{rel},\Gamma_2) \longrightarrow {\mathcal M}({\mathfrak Y}_1^\text{rel} \sqcup {\mathfrak Y}_2^\text{rel},\Gamma_1 \sqcup \Gamma_2)  \eqno(3.1) $$
to a closed substack of $M({\mathfrak W}_0,\Gamma)$, which glues two relative stable maps along their distinguished marked points. For each compatible pair $\eta = (\Gamma_1,\Gamma_2)$, there exist line bundles ${\mathbf L}_\eta$ on the total space ${\mathcal M}({\mathfrak W},\Gamma)$ with global sections ${\mathbf s}_\eta$ such that the vanishing locus of ${\mathbf s}_\eta$ is a closed substack ${\mathcal M}({\mathfrak W}_0,\eta)$ of ${\mathcal M}({\mathfrak W}_0,\Gamma)$ which is topologically identical to ${\mathcal M}({\mathfrak Y}_1^\text{rel} \sqcup {\mathfrak Y}_2^\text{rel},\Gamma_1 \sqcup \Gamma_2)$. 

Having briefly recalled the main objects, we can return to the specifics of the degeneration described in the previous section. For notational purposes, it is useful to introduce the usual space of stable maps $\overline{\mathcal M}_{0,0}(W,d)$ and the Chow variety $\text{Chow}_d(W)$. There are natural forgetful maps from ${\mathcal M}({\mathfrak W},\Gamma)$ to the former and from the former to the latter. The composition is denoted by $\gamma_{\mathfrak W}: {\mathcal M}({\mathfrak W},\Gamma) \to \text{Chow}_d(W)$.

There is a similar sequence of morphisms for each $Y_i$
$$ {\mathcal M}({\mathfrak Y}_i^{\text{rel}},\Gamma_i) \longrightarrow \prod_{v \in V(\Gamma_i)} \overline{\mathcal M}_{0}(Y_i,\beta(v)) \longrightarrow \prod_{v \in V(\Gamma_i)} \text{Chow}_{\beta(v)}(Y_i), \eqno(3.2) $$
where $v$ ranges over all vertices of the graph specified by $\Gamma_i$ and $\beta:V(\Gamma_i) \to \mathrm{H}_2(Y_i,{\mathbf Z})$ gives the class in $Y_i$ of the image of each component.

Let $\eta = (\Gamma_1,\Gamma_2)$ be a compatible topological type with $g_1 \equiv g_2 \equiv 0$. The expectation that ${\mathcal M}({\mathfrak W}_0,\eta) \cong_{\text{top}} {\mathcal M}({\mathfrak Y}_1^\text{rel} \sqcup {\mathfrak Y}_2^\text{rel},\Gamma_1 \sqcup \Gamma_2)$ has dimension zero is unrealistic because of the presence of multiple covers. However, we may ask whether its image under $\gamma_{\mathfrak W}$ has dimension zero -- this is false in general as well. Nevertheless, three questions arise naturally:

(Q1) Are the images of the closed fibers of $\mathbf{q}_i$ under $\gamma_{\mathfrak W}$ finite?

(Q2) Do the images of $\mathbf{q}_i$ have the expected dimensions?

(Q3) Do these images intersect dimensionally transversally for $i=1,2$?

\noindent  The answers to (Q1) and (Q2) are affirmative. Unfortunately, the answer to the third question is in general negative. However, the third point turns out to be of purely combinatorial nature and as such, it is possible to determine algorithmically when it does holds.

\subsection{Relative stable maps to ${\mathfrak Y}_1$ and ${\mathfrak Y}_2$} 

In this section, we want to describe the moduli spaces ${\mathcal M}({\mathfrak Y}_i^{\text{rel}},\Gamma_i)$ of relative stable maps to the ${\mathfrak Y}_i$, when $g_1 \equiv g_2 \equiv 0$, i.e. all connected components have arithmetic genus $0$. Since there are no nonconstant maps from a rational curve to a curve of geometric genus $1$, all irreducible components of any relative stable map are sent either inside the fibers of the maps 
$$ \varphi_i:Y_i \longrightarrow E_i $$
constructed in section 2, or to the lines in the rulings of the intermediary ruled varieties reviewed earlier, possibly as multiple covers. Formally, we have a morphism 
$$ \varphi^{\mathcal M}_i:{\mathcal M}({\mathfrak Y}_i^{\text{rel}},\Gamma_i) \longrightarrow E_i^{V(\Gamma_i)} $$
specifying the $\varphi_i$-fiber containing each connected component. We will show that the conditions the $r$-tuples of images in $S$ of the distinguishes marked points need to satisfy boil down to a collection of linear equations in $\text{Pic}(E)$. These conditions are obtained by restricting information about rational equivalence on the surfaces $\varphi_i^{-1}(p)$ to information about rational equivalence on the corresponding copy of $E_j \subset \varphi_i^{-1}(p)$.

First, we will address (Q1) above. To avoid mentioning the awkward maps to Chow varieties at each step, we introduce some ad hoc terminology. In this language, question (Q1) becomes: are the distinguished evaluation morphisms $\mathbf{q}_i$ cycle-finite? We prove that this is the case.

\begin{term}
Let $\pi:C \to S$ be a flat family of semistable curves over a ${\mathbf C}$-scheme (or perhaps DM-stack) $S$ and $f:C \to X \times_{\mathbf C} S$ a family of morphisms to some projective variety $X$, which induces a map $S \to \text{Chow}(X)$. Let $g:S \to Y$ be a proper morphism to a scheme $Y$. We say that $g$ is \textit{cycle-finite} relative to the family of maps $C \to X$ if the image in $\text{Chow}(X)$ of any closed fiber $g^{-1}(y)$ is zero-dimensional.
\end{term}
 
\begin{proposition}
The distinguished evaluation morphisms $\mathbf{q}_i$ are cycle-finite relative to the universal family over ${\mathcal M}({\mathfrak Y}_i^{\text{rel}},\Gamma_i)$ mapping into $Y_i$.
\end{proposition}

\begin{proof}
Roughly, positive-dimensional families would be forced to split off a component in $S$, which is impossible since $S$ doesn't contain rational curves. 

Denote the universal family by $\pi_1: {\mathcal M}({\mathfrak Y}_i^{\text{rel}},\Gamma'_i) \to {\mathcal M}({\mathfrak Y}_i^{\text{rel}},\Gamma_i)$. The notation $\Gamma'_i$ denotes a topological type identical to $\Gamma_i$, with the addition of a single ordinary marked point. Let 
$ \text{ev}_{1}: {\mathcal M}({\mathfrak Y}_i^{\text{rel}},\Gamma'_i) \to Y_i$ be the evaluation morphism associated to this point. Let $s = (s_1,...,s_r) \in S^r$ be a closed point. Our claim is essentially that the image $\Sigma$ of $\pi_1^{-1}({\mathbf q}_i^{-1}(s))$ under $\text{ev}_1$ in $Y_i$ has dimension one. The crucial observation is this: since $S = E \times E$ contains no rational curves, all components of the relative stable maps which map to the intermediary ruled threefolds over $S$ will be mapped as covers of some fibers of these rulings; this implies that $\Sigma \cap S$ is supported on $\{s_1,...,s_r\}$, so $\Sigma$ can't have dimension $2$ or more. 
\end{proof}

Choose once and for all an arbitrary hyperflex point of $E \subset {\mathbb P}^4$ to be the zero element of $(E,+)$. This gives convenient identifications of all connected components of $\text{Pic}(E)$ with $E$ and in particular, a preferred isomorphism $E_1 \cong E$. Consider the pushforward maps $\varphi_{i*} : \text{H}_2(Y_i,{\mathbf Z}) \to \text{H}_2(E_i,{\mathbf Z}) \cong {\mathbf Z}$. Let $K_1 \cong {\mathbf Z}[\text{fiber}] \oplus {\mathbf Z}[\mathrm{directrix}] \subset \text{H}_2(Y_1,{\mathbf Z})$ generated by the classes of the directrix respectively any line in the ruling of any fiber of $\varphi_1$ and $K_2$ generated by the classes of all lines on any smooth fiber of $\varphi_2$. Then $\mathrm{Im}(\beta_i) \subset K_i$.

For any vertex $v \in V(\Gamma_i)$, we denote by $\Gamma_i|_v$ the restriction of the topological data to $v$, i.e. the data consisting of the graph with the unique vertex $v$, the roots $R(v)$ which were previously attached to $v$, the class $\beta_i(v) \in \mathrm{H}_2(Y_i,{\mathbf Z})$, the genus function $g_i(v)=0$ and no legs. Consider the moduli stack ${\mathcal M}({\mathfrak Y}_i^\text{rel},\Gamma_i|_v)$ of relative stable maps of this topological type. There is a map 
$$ \varphi_{i|v}^{\mathcal M}: {\mathcal M}({\mathfrak Y}_i^\text{rel},\Gamma_i|_v) \longrightarrow E_i. $$
Let $\smash{ \Gamma'_i|_v} $ denote the topological type which is identical to $\Gamma_i|_v$, with the addition of one ordinary marked point. We obtain the universal family  $\smash{ {\mathcal M}({\mathfrak Y}_i^\text{rel},\Gamma'_i|_v) }$ over the previous moduli space. Recall the map $ {\mathcal M}({\mathfrak Y}_i^\text{rel},\Gamma'_i|_v) \to {\mathcal M}({\mathfrak Y}_i^\text{rel},\Gamma_i|_v) \times_{E_i} Y_i$ specifying the fibers of the components. Define the map
$$ \omega_i: {\mathcal M}({\mathfrak Y}_i^\text{rel},\Gamma_i|_v) \longrightarrow {\text{Pic}}^{\delta(v)} (E_j) \eqno(3.3) $$ 
$$ \left( C,f,(q_\alpha )_{\alpha \in R(v)} \right) \mapsto \sum_{\alpha \in R(v)}\mu(\alpha)f(q_\alpha) \in {\text{Pic}}^{\delta(v)} (E_j), $$
where $\delta(v) = \sum_{\alpha \in R(v)} \mu(\alpha)$, i.e. the divisor class cut out by $f(C)$ on $\smash{ \{ \varphi_{i|v}^{\mathcal M}(f) \} \times E_j }$. The definition in families is similar, just harder to type. The first step is to understand the image of the map $\smash{ \omega_i \times \varphi_{i|v}^{\mathcal M} }$ in $\smash{ \text{Pic}^{\delta (v)}(E_j) \times E_i }$. Identify the target with $E_j \times E_i = E_1 \times E_2$ and let $(x,y)$ be coordinates on $E_1 \times E_2$.

\begin{lemma} 
$(a)$ Let $i=1$ and $\beta_1(v)=k_f[\mathrm{line}] + k_\infty[\mathrm{directrix}] \in K_1$. Then the image is contained inside the locus in $E_1 \times E_2$ given by the equation $y - (k_f-2k_\infty)x = 0$.

(b) Let $i=2$ and $\beta_2(v) \in K_2$ of degree $k_l$ relative to ${\mathscr O}_{X_{0,2}}(1)$. Then the image is contained in the union of the curve of equation $3x = k_ly$ with $E_1 \times \{p_1,...,p_{25}\}$.
\end{lemma}

\begin{proof}
$(a)$ This follows from part (a) of Lemma 2.4. For any $D \in E_1$, the divisor class ${\mathscr O}_\Sigma(k_f f + k_\infty e_\infty)$ on the fiber $\Sigma = \varphi_1^{-1}(D)$ restricts on the corresponding copy of $E_2 = E$ sitting inside $\Sigma$ to ${\mathscr O}_{{\mathbb P}^4}(k_\infty)|_E \otimes {\mathscr O}((k_f-2k_\infty)D)$. Given that the implicit isomorphism of the corresponding connected component of $\mathrm{Pic}(E)$ with $E$ is induced by a hyperflex point, the first terms disappears, so the image lies inside the locus $y = (k_f-2k_\infty)x$.

$(b)$ Similarly, this follows from part (b) of Lemma 2.4 and the discussion in \S2.3. Let $p \in E_2=E$. Assume that $p \notin \{p_1,p_2,...,p_{25}\}$. Then the fiber $\Sigma = \varphi_2^{-1}(p)$ of $\varphi_2$ is a smooth cubic surface, as in \S2.3. By the discussion in that section, any divisor class on $\Sigma$ which pushes forward to $k_l[\text{line}]$ in $\text{H}_2(Y_2,{\mathbf Z})$ has the property that its cube is the class ${\mathscr O}_{{\mathbb P}^4}(k_l)|_E \otimes {\mathscr O}_E(-2k_lp)$, so, with the implicit identifications, we get the desired constrain $3(x-k_ly)+2k_ly = 0$. The reason for the term $x-k_ly$ is the following: the coordinate $x$ on $E_1 \times E_2$ is actually the coordinate $\lambda$ on the canonical $E \times E$ and the corrections will stack up in $\text{Pic}^{k_l}(E_1)$.
\end{proof}

Let $L(\Gamma_{i|v})$ be the image of $\omega_i \times \varphi^{\mathcal M}_{i|v}$ in $\mathrm{Pic}^{\delta(v)} {(E_j)} \times E_i$. Let $\mathbf{q}_{i|v}$ be the evaluation morphism restricted to $v$ and consider the map which associates to a collection of closed points on $E_j$ indexed by $R(v)$, the linear equivalence class of the corresponding $\mu$-weighted divisor on $E_j$, i.e. 
$$ h:E_j^{R(v)} \longrightarrow \mathrm{Pic}^{\delta(v)}(E_j) $$
$$ (y_\alpha)_{\alpha \in R(v)} \mapsto \sum_{\alpha \in R(v)} \mu(\alpha) y_\alpha \in \text{Pic}^{\delta(v)} {(E_j)}. $$
Motivated by the restrictions imposed by Lemma 3.3 on the configurations of distinguished marked points, define the following loci in $\smash{ (E_1 \times E_2)^{R(v)} }$:
$$ A(\Gamma_{i|v}) = (\text{Id}_{E_j^{R(v)}} \times \text{Diag}_{E_i})\left( (h \times \text{Id}_{E_i})^{-1}L(\Gamma_{i|v}) \right), $$
for $i=1,2$. Note that $\smash{ A(\Gamma_{i|v}) }$ is an equidimensional union of abelian subvarieties of the abelian variety $\smash{ (E_1 \times E_2)^{R(v)} }$, of dimension $ \smash{ \dim ( E_j^{R(v)} \times E_i ) - 1 = |R(v)|}$, and that the image of ${\mathbf q}_{i|v}$ is contained in $A(\Gamma_{i|v})$. Let $A(\Gamma_{i}) \subset (E_1\times E_2)^r$ be the direct product of all $A(\Gamma_{i|v})$ over all $v \in V(\Gamma_i)$. Then the image of ${\mathbf q}_i$ is contained in $A(\Gamma_{i})$, which is a union of abelian subvarities of $(E_1 \times E_2)^r$ of pure dimension $\smash { \sum_{v \in V(\Gamma_i)} |R(v)| = r }$. Note that $A(\Gamma_{1})$ and $A(\Gamma_{2})$ have complementary dimensions in $S^r$, as we should expect.
\begin{center}
\begin{tikzpicture}
\matrix [column sep  = 20mm, row sep = 15mm] {

	\node (sw) {$(E_1 \times E_2)^{R(v)}$}; &
	\node (s) {$E_j^{R(v)} \times E_i$}; &
	\node (se) {$\text{Pic}^{\delta(v)} {(E_j)} \times E_i$}; & \\	
	
	& \node (n) {${\mathcal M}({\mathfrak Y}_i^\text{rel},\Gamma_i|_v)$}; \\
};
\draw[-{Stealth[length=2mm]}, thin] (n) -- (sw);
\draw[-{Stealth[length=2mm]}, thin] (n) -- (s);
\draw[-{Stealth[length=2mm]}, thin] (n) -- (se);
\draw[-{Stealth[length=2mm]}, thin] (s) -- (sw);
\draw[-{Stealth[length=2mm]}, thin] (s) -- (se);

\node at (-3.7,-0.2) {$\mathbf{q}_{i|v}$};
\node at (-3,1.5) {$\text{Id}_{E_j^{R(v)}} \times \text{Diag}_{E_i}$};
\node at (1,1.5) {$h \times \text{Id}_{E_i}$};
\node at (1.9,-0.2) {$\omega_i \times \varphi^{\mathcal M}_{i|v}$};

\end{tikzpicture}
\end{center}

We conclude this subsection with a remark -- or rather, a question for the reader. Given any abelian variety and any two abelian subvarieties of complementary dimensions, either the two abelian subvarieties intersect transversally, or the intersection of the associated homology classes is zero. This can be seen, for instance, from the excess intersection formula. 

From the point of view of the degeneration formula \cite{[Li02]}, this seems to have intriguing consequences. Any hypothetical family of limits which is not cycle-finite is necessarily associated with the vanishing of the corresponding part of the contribution to the virtual count. However, a huge technical difficulty is the following: for $i=2$ the image consists not of an abelian subvariety, but a singular union of subvarieties and it is entirely unclear whether it is possible to separate the contributions from each individual component in a geometrically meaningful way.  

\subsection{Combinatorics of the degenerate maps}
In this section, we write explicitly the combinatorial laws governing the stable maps to ${\mathfrak W}_0$. Let $\Gamma_1$ and $\Gamma_2$ be compatible topological types, i.e. $\Gamma_1$ and $\Gamma_2$ have the same number $r$ of distinguished marked points and the identically indexed distinguished marked points have the same $\mu$-weight. As we said, we will only be concerned with the case $g_1 \equiv g_2 \equiv 0$, which covers all cases if $g=0$ and a fraction of the cases when $g>0$.    

To avoid notation becoming too heavy, we will work with only one map instead of families. Fix a single stable map $(C,f,q_1,...,q_r)$ obtained by gluing a relative stable map in ${\mathcal M}({\mathfrak Y}_1^\mathrm{rel},\Gamma_1)$ and a relative stable map in ${\mathcal M}({\mathfrak Y}_2^\mathrm{rel},\Gamma_2)$. Let $G$ be the graph with vertices $V(G) = V(\Gamma_1) \cup V(\Gamma_2)$ and $r$ edges obtained by gluing the roots with identical indexing. Then $\text{rank } \mathrm{H}^1(G) = g$, since 
$$ g = \text{rank } \mathrm{H}^1(G)+\sum_{v \in V(\Gamma_1)} g_1(v)+ \sum_{v \in V(\Gamma_2)} g_2(v) $$
and we are assuming $g_1 \equiv g_2 \equiv 0$. From now on, we will call the elements of $V(\Gamma_1)$ red vertices and the elements of $V(\Gamma_2)$ blue vertices. As in Lemma 3.3, let $k_f,k_\infty:V(\Gamma_1) \to {\mathbf N}$ such that $\beta_1 \equiv k_f[\text{line}] + k_\infty[\text{directrix}]$ and $k_l:V(\Gamma_2) \to {\mathbf N}$ such that $\beta_2$ has degree $k_l$ relative to ${\mathscr O}_{X_{0,2}}(1)$. 

Let $\lambda:V(G) \to (E,+) \cong ({\mathbf R}^2/{\mathbf Z}^2,+)$ be the map giving the components of $\varphi_1^{\mathcal M}([f])$ and $\varphi_2^{\mathcal M}([f])$. It is possible to merge $i=1$ and $i=2$ as we are simultaneously identifying $E_1$ and $E_2$ with $E$. The point is that Lemma 3.3 and the discussion thereafter imposes $|V(G)|$ linear conditions on $\lambda$, which are usually linearly independent. Denote the set of neighbors of $v$ in $G$ by $V(v)$. 

First, consider the case when $v$ is a red vertex. Then
$$ \sum_{w \in V(v)} \mu([vw]) \lambda(w) = \left( k_f(v) - 2k_\infty(v) \right) \lambda(v), $$
with the numerical constrain on the multiplicities
$$ \sum_{w \in V(v)} \mu [vw]  =  2k_f(v) + k_\infty(v),  $$
since $(k_f[\text{line}] + k_\infty[\text{directrix}] \cdot -K_{{\mathbb F}^1})_{{\mathbb F}_1} = 2k_1^f(v) + k_1^\infty(v)$ by a straightforward calculation.  

Second, consider the case when $v$ is blue. Then either $\lambda(v) \in \{p_1,p_2,...,p_{25}\}$, when nothing more can be said, or otherwise, again from the computation in Lemma 3.3 and the following discussion, we get
$$ 3 \sum_{w \in V(v)} \mu ([vw]) \lambda(w) = k_l(v) \lambda(v)  $$
and the numerical constrain on the weights
$$ \sum_{w\in V(v)} \mu ([vw]) = k_l(v). $$
Note that the latter constrain still holds in the case $\lambda(v) \in \{p_1,p_2,...,p_{25}\}$. We introduce some (fairly standard) notation to state the equations above in a slightly more pleasant form. 

\begin{definition}
Let $G$ be any connected graph without isolated vertices and let $\mu:E(G) \to {\mathbf Z}^{+}$ be a weight function on the set of edges. First, we naturally define the weighted degrees of the vertices as
$$ \deg_\mu(v) = \sum_{w\in V(v)} \mu ([vw]). $$
For any function $f$ defined on the set of vertices $V(G)$ with values in some abelian group, we define the $\mu$-weighted unnormalized $G$-Laplacian of $f$ by
$$ \Delta^u_\mu f(v) := \sum_{w \in V(v)} \mu([vw])(f(v)-f(w)) $$
The normalized $\mu$-weighted $G$-Laplacian $\Delta_\mu$ is defined by the same formula, but dividing by $\deg_\mu(v)$, assuming that $\deg_\mu(v)$ is invertible in the target group. If $\mu \equiv 1$, we suppress this subscript.
\end{definition}

\noindent A basic computation allows us to rewrite the conditions above as follows:
$$ \begin{cases}
\Delta_\mu^u\lambda(v) = \left( 3\deg_\mu(v) -5k_f(v) \right)\lambda(v)  & \text{if } v \text{ is red,} \\[2ex]
3 \Delta_\mu^u\lambda(v) = 2\deg_\mu(v) \lambda(v) \text { or } \lambda(v) \in \{p_1,p_2,...,p_{25}\} &\text{if } v \text{ is blue.}
\end{cases} \eqno(3.4) $$
In conclusion, (cycle-)finiteness of ${\mathcal M}({\mathfrak W}_0,\eta)$ boils down to the system of linear equations in ${\mathbf R}^2/{\mathbf Z}^2$, or equivalently ${\mathbf R}/{\mathbf Z}$, having only finitely many solutions. For $g=0$, the smallest counterexample to finiteness occurs in degree $d=6$; see Example 3.6 below.

Having completed the analysis of the degenerate stable maps, we would like to briefly return to the original question, that of finding the limits in the classical sense. First, we show two examples to illustrate the following facts: (1) different degenerate stable maps mapping to the same $1$-cycle on $W_0$ may behave differently; and, more importantly (2) some degenerate stable maps do not correspond to classical limits at all. Then, we propose a speculation on the position of the distinguished marked points for classical limits. 

\begin{example}
For $d=1$, $g=0$, consider a graph $G$ with a red vertex $R$ and a blue vertex $B$ such that $k_f(R) = 1$, $k_\infty(R) = 0$, $k_l(B) = 2$ and weight $2$ on the edge. 

Geometrically, some of these maps can be described as follows: the red component maps one-to-one to a line in $Y_1$ which is tangent to $E \times E$, while the blue component maps two-to-one to a line in $Y_2$ and the distinguished marked point is a ramification point. By moving the remaining ramification point, we obtain a one-parameter family of such maps. Of course, for dimension reasons, only finitely many can be real limits.
\end{example}

\begin{example} (This example was found by G. Bujokas.) For $r=2$, we may find examples in which the system (3.4) is special. Let $d=6,g=0$ and $G$ have two red vertices $R_1,R_2$ and a blue vertex $B$ of "smooth" type. Let $k_f(R_i) = 2$, $k_\infty(R_i) = 1$, $k_l(B) = 10$ and let both edges have weight $5$. A rough picture is provided below; the point $b = \lambda(B)$ denotes the point of $E_2$ whose fiber contains the blue component and similarly for the points $r_1=\lambda(R_1),r_2=\lambda(R_2) \in E_1$. 

An immediate calculation shows that (3.4) is special in this case. Geometrically, the red vertices correspond to twisted cubics which are sections of the corresponding ${\mathbb F}_1$ surfaces by the same osculating hyperplane in ${\mathbb P}^4$ to a hyperflex of $E$.
\end{example}

Let $\Lambda$ be the countable subset of $E$ consisting of points which are rational combinations of the $25$ points $p_1,...,p_{25}$, that is, $kq \sim l_1 p_1 + ... + l_{25} p_{25}$ for $k,l_1,...,l_{25} \in {\mathbb Z}$ and $k \neq 0$. Recall the map $\zeta$ introduced in section 1.3. Motivated by the comment at the end of section 3.2 and the form of the system (3.4), we propose the following.

\begin{question}
Let $q$ be any intersection point of some limit of rational curves with $E \times E$. Then, is it necessarily true that $q \in \Lambda \times \Lambda$?
\end{question}

Of course, if the answer is affirmative, the finiteness part of the Clemens conjecture follows. Further evidence in favor of Question 3.7 comes from the following interesting observation of X. Chen: $q$ necessarily belongs to the union of the (countably many) smooth genus $1$ curves on $E \times E$ which pass through points of $\Lambda \times \Lambda$. However, this is actually a feature of arbitrary degenerate genus $0$ stable maps and it can be proved by a clever combinatorial argument involving only (3.4).

\begin{center}
\begin{tikzpicture}[scale=0.9, every node/.style={scale=0.9}]

\def\xa{1}
\def\xb{3}
\def\ya{1}
\def\yb{0}
\def\k{1}
\def\j{3.5}
\def\l{3}

\draw [dotted] (\xa + \xb,\ya + \yb) -- (\xa-\xb,\ya-\yb);

\draw (\xa-\xb+0.5,\ya-\yb) -- (\xa-\xb,\ya-\yb);
\draw (\xa-\xb+2.5,\ya-\yb) -- (\xa-\xb+3,\ya-\yb);
\draw (\xa-\xb+5,\ya-\yb) -- (\xa-\xb+6,\ya-\yb);

\draw (\xa + \xb,\ya + \yb) -- (-\xa+\xb,-\ya+\yb);
\draw (-\xa - \xb,-\ya - \yb) -- (\xa-\xb,\ya-\yb);
\draw (-\xa - \xb,-\ya - \yb) -- (-\xa+\xb,-\ya+\yb);

\draw (\xa + \xb-\k,\ya + \yb) -- (-\xa+\xb-\k,-\ya+\yb);
\draw (\xa + \xb-\k,\ya + \yb +\l) -- (-\xa+\xb-\k,-\ya+\yb +\l );
\draw (\xa + \xb-\k,\ya + \yb) -- (\xa + \xb-\k,\ya + \yb +\l);
\draw (-\xa+\xb-\k,-\ya+\yb) -- (-\xa+\xb-\k,-\ya+\yb +\l );

\draw (\xa + \xb-\j,\ya + \yb) -- (-\xa+\xb-\j,-\ya+\yb);
\draw (\xa + \xb-\j,\ya + \yb +\l) -- (-\xa+\xb-\j,-\ya+\yb +\l );
\draw (\xa + \xb-\j,\ya + \yb) -- (\xa + \xb-\j,\ya + \yb +\l);
\draw (-\xa+\xb-\j,-\ya+\yb) -- (-\xa+\xb-\j,-\ya+\yb +\l );

\draw [dotted] (-\xb,0) -- (\xb,0);

\draw (-\xb,0) -- (-\xb+1.5,0);
\draw (-\xb+2.5,0) -- (-\xb+4,0);
\draw (-\xb+5,0) -- (-\xb+6,0);

\draw [dotted] (-\xb,0) -- (-\xb,-\l);

\draw (-\xb,-1) -- (-\xb,-\l);
\draw (\xb,0) -- (\xb,-\l);
\draw (-\xb,-\l) -- (\xb,-\l);

\draw [very thick, blue] (-0.5,0) to [in= 90, out = 180]++(-1,-2);
\draw [very thick, blue] (-0.5,0) to [in= 180, out = 0]++(1,-2);
\draw [very thick, blue] (2,0) to [in= 90, out = 0]++(0.7,-2);
\draw [very thick, blue] (2,0) to [in= 0, out = 180]++(-1.5,-2);

\draw [very thick, red] (-0.5,0) to [in= 270, out = 45]++(0.7,3);
\draw [very thick, red] (-0.5,0) to [in= 270, out = 225] ++(-0.7,2);

\draw [very thick, red] (2,0) to [in= 270, out = 45]++(0.7,3);
\draw [very thick, red] (2,0) to [in= 270, out = 225] ++(-0.7,2);

\node at (-0.5,0) {$\bullet$};
\node at (2,0) {$\bullet$};

\node at (-2.8,1.8) {$\varphi_1^{-1}(r_2) \subset Y_1$};
\node at (4.2,1.8) {$\varphi_1^{-1}(r_1) \subset Y_1$};

\node at (-0.6,0.6) {$\mu=5$};
\node at (1.9,0.6) {$\mu=5$};

\node at (2.4,2) {$R_1$};
\node at (-0.1,2) {$R_2$};
\node at (1.2,-2) {$B$};
\node at (-1.7,-2.5) {$\varphi^{-1}_2(b) \subset Y_2$};

\node at (-3,-0.6) {$S$};

\end{tikzpicture}
\end{center}

\begin{center}
\begin{tiny}
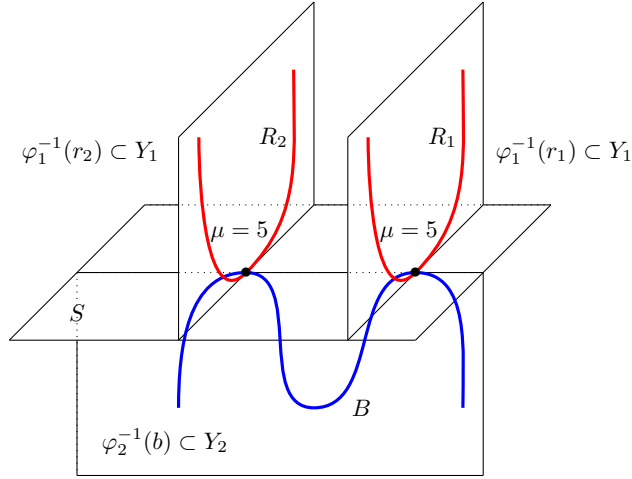
\textbf{Fig. B.} The picture for Example 3.6.
\end{tiny}
\end{center}

\subsection{Existence of rigid stable maps}
In this section, we prove Theorem 1.1 by explicitly exhibiting rigid degenerate stable maps to ${\mathfrak W}_0$, of arithmetic genus $g$ and degree $d \gg_g 0$ with certain additional properties and invoking the existence of a perfect obstruction theory on ${\mathcal M}({\mathfrak W},\Gamma)$ to infer that they are indeed limits of rigid stable maps to nearby fibers of $W \to {\mathbb A}^1$. Smoothness of the source follows easily from the smoothness of the connected components of the maps to ${\mathfrak Y}_1^\text{rel}$ and ${\mathfrak Y}_2^\text{rel}$. 

The amount of freedom in choosing the topological types is remarkable and reassuring. The first main assumption is $\beta_1 \equiv [\text{line}]$, that is, we are assuming that all red vertices correspond to fibers of the rulings in the corresponding Hirzebruch surfaces. Furthermore, we assume $\mu \equiv 1$, hence all red vertices have degree $2$. Then the first branch of equation $(3.4)$ simply reads $\lambda(v) = \lambda(w_1) + \lambda(w_2)$, where $v$ is any red vertex and $w_1,w_2$ are its two neighbors.

Let $G'$ be the graph obtained from $G$ by suppressing all red vertices and replacing any length $2$ chain in $G$ between two blue vertices with an edge. Then the second branch of equation (3.4) becomes simply $3\Delta^u \lambda(v) = 5\deg (v) \lambda(v)$, where the discrete Laplacian is now taken relative to $G'$. We further assume that $\deg v \leq 3$ for all vertices of $G'$ and proceed to construct the degenerate curve.

\begin{claim}
For $d,g \geq 1$, $d \gg_g 0$, there exists a set of data as follows:

$\bullet$ A connected graph $G=(V,E)$ with $|E| = d$ and $h^1(G)=g$, hence $|V| = d-g+1$. We require that $\deg v \leq 3$ for all $v \in V$.

$\bullet$ If ${\mathbf K}$ is some field of characteristic $\neq 2,3$, consider the discrete equation
$$ \Delta \lambda = \frac{5}{3}\lambda \eqno(3.5) $$
in $\text{Fun}(V,{\mathbf K})$, the ${\mathbf K}$-vector space of functions from $V$ into ${\mathbf K}$. We require that $(3.5)$ has only the trivial solution $\lambda = 0$ when ${\mathbf K} = \overline{\mathbf Q}$, i.e. $5/3$ is not an eigenvalue of the Laplacian.

$\bullet$ A prime number $p \geq 7$ and a solution $\hat{\lambda}$ of $(3.5)$ for ${\mathbf K} = {\mathbf F}_p$ such that:

1. $\hat{\lambda}(v) \neq \hat{\lambda}(w)$ if $1 \leq \text{dist}(v,w) \leq 2$; and

2. $\hat{\lambda}$ is a \textit{strongly irreducible} solution of $(3.5)$, in the following sense. For any induced subgraph $H \subset G$, consider the similar discrete equation
$$ F_H \lambda :=  \Delta_H \lambda - \frac{5}{3}\lambda = 0. $$
Then we require that $\hat{\lambda}$ has the following property: if $F_H \hat{\lambda}(v) = 0$, then $\deg_H v$ is equal to either $\deg_G v$ or $0$.
\end{claim}

\begin{lemma}
If no two degree $3$ vertices of $G$ are adjacent, then $(3.5)$ only has the trivial solution $0$ for ${\mathbf K} = {\mathbf Q}$, or equivalently, $\overline{\mathbf Q}$.
\end{lemma}

\begin{proof}
Let ${\mathbf K} = {\mathbf Q}$ and let $v_3:{\mathbf Q} \to {\mathbf Z} \cup \{+\infty\}$ be the valuation at the prime $3$. Pick any solution $\lambda$ of $(3.5)$. Then we have
$$ v_3(\lambda(v)) = v_3\left(\sum_{w \in V(v)} \lambda(w) \right) - v_3(\deg v) + 1 \geq \min_{w \in V(v)} v_3(\lambda(w)) + 1 - v_3(\deg v). $$
Hence any vertex $v$ has a neighbor with strictly smaller $v_3 \circ \lambda$, if $\deg v \leq 2$ and $\lambda(v) \neq 0$, respectively no greater $v_3 \circ \lambda$, if $\deg v = 3$. This easily implies $v_3 \circ \lambda \equiv +\infty$, so $\lambda \equiv 0$.
\end{proof}

\begin{claim}
There exists a prime $p \geq 7$ and three graphs $G_1,G_2,G_3$ with eigenvectors $\hat{\lambda}_i$ in ${\mathbf F}_p$ satisfying the conditions in Claim 3.8 and Lemma 3.9 ($d$ and $g$ are not fixed) and two elements $a,b \in {\mathbf F}_p$ such that:

$\bullet$  $|E(G_1)|$ and $|E(G_2)|$ are coprime;

$\bullet$ $h^1(G_1) = h^1(G_2) = 1$ and $h^1(G_3) = 2$;

$\bullet$ all graphs $G_i$ contain an edge $[vw]$ such that $v$ and $w$ both have degree $2$ in $G_i$ and $\hat{\lambda}_i(v) = a$, $\hat{\lambda}_i(w)=b$.
\end{claim}

\begin{proof}
The proof is by explicit example found\footnote{A much simpler computer program for the purpose of verifying these examples is available on the author's website \url{https://math.ucdavis.edu/~azahariuc/research/main[1].cpp}. This verification suffices logically to complete the proof of Theorem 1.1.} by computer, but verifiable by hand. However, the point is that, while such examples can't be particularly simple or enlightening, there is no conceivable reason why they wouldn't exist and, in fact, there is a very wide variety of such examples. 

The prime number is $p=23$. The graph $G_1$ is a length $11$ cycle. The eigenvector $\lambda:V(G_1) \to {\mathbf F}_{23}$ of $\triangle$ with eigenvalue $5/3$ is simply a way to associate mod $23$ integers to the vertices. For $G_1$, the entries are, in order, $1, 2, 4, 8, 16, 9, 18, 13, 3, 6, 12$.

The graph $G_2$ has $18$ edges and vertices, and consists of a length $15$ cycle with $3$ "spikes"; the entries are $1,2,4,8(11),15,18,7,11,9(21),19,4,14,8,$ $6(14),12$.

Finally, $G_3$ has $43$ edges and is somewhat more complicated. There are two special vertices, both of degree $3$, with entries $9$ and $7$. Then, there are three "bridges" between these two vertices, as follows. The first bridge consists of a chain with $20$ vertices excluding the ends and $4$ spikes. The entries of the eigenvector are in order $15, 17, 16 (22), 21, 2, 7, 4, 3 (7), 6, 12, 1, 2, 4 (17),$ $ 19, 9, 15, 17, 16 (22), 21, 2$. The sole purpose of the long repeating sequence is to satisfy the third bullet in the claim. The second bridge is $16, 8 (11), 3, 11,$ $13, 10 (8), 5, 14$ and the third bridge is $20, 18, 2, 10 (8), 16$. \end{proof}

\begin{proof}[Proof of Claim 3.8] We will use the graphs $G_1,G_2,G_3$ constructed in Claim 3.10 as the building blocks for constructing $G$. For each $G_i$, we may define a graph $H_i$ with $|E(H_i)| = |E(G_i)|+1$ and $h^1(H_i) = h^1(G_i)-1$ obtained by duplicating the edge $[vw]$ and then separating the two ends, i.e. $H_i$ will begin and end with a copy $[vw]$. 

Starting with some $G_i$, we may insert copies of various $H_j$'s at the edges which are copies of $[vw]$. Note that we may glue the functions from vertices to ${\mathbf F}_p$ and crucially, the property of being eigenvectors with eigenvalue $5/3$ is preserved since any collections of vertices consisting of some arbitrary vertex and all its neighbors is contained isomorphically in some $G_i$. It is also not hard to see that the conditions $\hat{\lambda}(v) \neq \hat{\lambda}(w)$ if $1 \leq \text{dist}(v,w) \leq 2$ respectively $\deg v \neq 3$ if $v$ has some neighbor of degree $3$ are preserved as well through this process. To ensure $h^1(G) =g$, it is enough to insert $H_3$ precisely $g-1$ times throughout the process. Finally, since $|E(G_1)|$ and $|E(G_2)|$ are coprime, the elementary coin problem ensures that $G$ may have any sufficiently high number of edges.
\end{proof}

\begin{lemma}
For $d \gg g$, $g \geq 1$ there exists a stable map $f_0:C_0 \to W_0$, $[f_0] \in {\mathcal M}({\mathfrak W}_0,\Gamma)$ obtained by gluing two relative stable maps $f_1 \in {\mathcal M}({\mathfrak Y}_1^\text{rel},\Gamma_1)$ and $f_2 \in {\mathcal M}({\mathfrak Y}_2^\text{rel},\Gamma_2)$ mapping to $Y_1=Y_1[0]$ (meaning that no expansions of $(Y_1,S)$ are required, cf. \cite{[Li02]}), respectively $Y_2=Y_2[0]$, such that all connected components of the source of $f_i$ are smooth and are mapped generically one-to-one onto their images and $[f_0]$ is an isolated point of ${\mathcal M}({\mathfrak W}_0,\Gamma)$. 
\end{lemma}

\begin{proof} We will avoid going into the details of this proof, since the desired degenerate map is obtained simply by reversing the discussion so far, with few significant new ingredients. First, we use the graphs constructed in Claim 3.8 as the $G'$ in the discussion at the beginning of this section. Using these, we may reconstruct $G$, the dual graph of $C$. Second, we use the solution $\smash{\hat{\lambda}}$ of (3.5) in ${\mathbf F}_p$ to reconstruct the desired solution $\lambda$ of (3.4) as follows: first, we extend $\smash{\hat{\lambda}}$ to $(E({\mathbf C}),+)$ by first embedding ${\mathbf Z}/p{\mathbf Z}$ diagonally into ${\mathbf Z}/p{\mathbf Z} \times {\mathbf Z}/p{\mathbf Z}$, which naturally sits inside $({\mathbf R}^2/{\mathbf Z}^2,+) \cong (E({\mathbf C}),+)$, then define $\lambda$ on the red vertices simply by $\lambda(v) = \lambda(w_1) + \lambda(w_2)$, where $w_1,w_2$ are the two blue neighbors of some red $v$.

Next, we construct the components $f_v:C_v \to Y_i$ of the stable map. We start with $v$ blue. Consider the divisor $D=\smash{\sum_{w \in V(v)} \lambda(w)}$ on the copy of $E$ sitting inside the cubic surface $\smash{\Sigma = \varphi_2^{-1}(\lambda(v))}$. The fact that $\lambda$ is a solution of (3.4) implies that we are in the situation of Lemma 2.11. Let $D_\Sigma$ be the divisor on $\Sigma$ such that $D_\Sigma \cap E = D$, as constructed in the said lemma. The crucial point (which will only be sketched) is that $D_\Sigma$ is irreducible. This is where the (combinatorial) strong irreducibility condition comes into play. Indeed, if the divisor were reducible, then either component would intersect the boundary $E$ at points which satisfy the analogous linear conditions in $(E({\mathbf C}),+)$ -- precisely what is ruled out by the strong irreducibility assumption. We are implicitly relying on the equivalence of the linear conditions for $\smash{\hat{\lambda}}$ and for $\lambda$. In conclusion, we can define $f_v:C_v \cong {\mathbb P}^1 \to Y_2$ to be the normalization of the divisor above. Finally, the distinguished marked points on $C_v$ are the preimages of $E$ under $f_v$. The fact that they are distinct is ensured by the condition $\hat{\lambda}(v) \neq \hat{\lambda}(w)$ if $\text{dist}(v,w) = 2$ in $G'$.

The case $v$ red is similar but less laborious. Again, the fact that $\lambda$ is a solution of $(3.4)$ ensures that we may find such a line. The fact that the distinguished marked points are distinct is ensured by $\hat{\lambda}(v) \neq \hat{\lambda}(w)$ if $v$ and $w$ are neighbors in $G'$. 

All in all, we may glue $f_1$ and $f_2$ to obtain a morphism $f_0:C_0 \to W_0$, where $C_0$ is obtained by gluing nodally all components of $C_v$ along the pair of identically indexed distinguished marked points. The fact that the distinguished marked points were obtained distinct directly from the construction above means that the target of $f_0$ is indeed $W_0$, i.e. no expansions are necessary. It is clear that $\text{Aut}(f_0)$ is trivial, so $f_0$ is stable. The crucial condition that $[f_0]$ is isolated in ${\mathcal M}({\mathfrak W}_0,\Gamma)$ follows from the second bullet in Claim 3.8.
\end{proof}

\begin{proof}[Proof of Theorem 1.1] If $g=0$, the statement is known by work of Katz, as explained in the introduction. If $g \geq 1$, the previous lemma all but completes the proof of the theorem. Since ${\mathcal M}({\mathfrak W},\Gamma)$ admits a perfect obstruction theory of the expected dimension \cite{[Li02]}, the irreducible component of ${\mathcal M}({\mathfrak W},\Gamma)$ containing the point $[f_0]$ must have dimension at least one. In fact, since the intersection of this locus with the central fiber has an isolated point at $[f_0]$, it must be of dimension precisely $1$ and, moreover, it must not be contained in the central fiber. It follows that we can find an ${\mathbb A}^1$-morphism $h:N \to {\mathcal M}({\mathfrak W},\Gamma)$ from a smooth affine curve $N$ with a distinguished point $x_0 \in N$ admitting a map $\psi:N \to {\mathbb A}^1$ such that $h(x_0) = [f_0]$ lies in the central fiber, but no other $h(x)$, $x \in N$ does. We may choose $N$ such that all $h(x_t)$ are isolated points of $\overline{\mathcal M}_{0,0}(W_{\psi(t)},d)$. 

Finally, we have to argue that the source $C_t$ of $h(x_t) = [(C_t,f_t)]$ is smooth for all $t$ in a punctured neighborhood of $x_0 \in N$, where $C \to N$ is the pullback to $N$ of the universal family of stable maps.
\begin{center}
\begin{tikzpicture}
\matrix [column sep  = 17mm, row sep = 8mm] {
	\node (nw) {$C$}; &
	\node (ne) {$W$};  \\
	\node (sw) {$N$}; &
	\node (se) {${\mathbb A}^1$}; \\
};
\draw[-{Stealth[length=2mm]}, thin] (nw) -- (ne);
\draw[-{Stealth[length=2mm]}, thin] (sw) -- (se);
\draw[-{Stealth[length=2mm]}, thin] (nw) -- (sw);
\draw[-{Stealth[length=2mm]}, thin] (ne) -- (se);

\node at (0,0.9) {$f$};
\node at (0,-0.4) {$\psi$};
\end{tikzpicture}
\end{center} 
The commutative diagram above induces a map $C \to W \times_{{\mathbb A}^1} N$. Any singular point of $C_0$ is a distinguished marked point $q$, so, in particular it is a node. The versal deformation space of a node is the germ at the origin of the family 
$$ \mathrm{Spec}({\mathbf C}[x,y,\tau]/(\tau-xy)) \to {\mathbb A}^1_\tau, $$ so it suffices to prove that $C$ is locally irreducible near $q$ in the complex-analytic topology. If this wasn't the case, then it would have two components $Z_1$ and $Z_2$, indexed such that $Z_i$ contains the branch of $f_0$ near $q$ mapping to $Y_i$. Then the pullback of the line bundle ${\mathscr O}_W(Y_1)$ with the section $1$ restricts on $Z_2$ to a line bundle with a section vanishing only at $q$, which is impossible. 
\end{proof}


\begin{thebibliography}{PP}

\bibitem{[Cl83]}
H. Clemens, \emph{Homological equivalence, modulo algebraic equivalence, is not finitely generated}, Pub. Math. IHES \textbf{58}, no. 1, 19--38 (1983)
\bibitem{[Cl86]}
H. Clemens, \emph{Curves on higher-dimensional complex projective manifolds}, Proc. Int. Cong. Math., Berkeley 634--640 (1986)
\bibitem{[Co05]}
E. Cotterill, \emph{Rational curves of degree 10 on a general quintic threefold}, Comm. Alg. \textbf{33}, no. 6, 1833--1872 (20).
\bibitem{[Co12]}
E. Cotterill, \emph{Rational curves of degree 11 on a general quintic 3-fold}, Quart. J. Math. \textbf{63}, no. 3, 539--568 (2012)
\bibitem{[EH87]}
D. Eisenbud and J. Harris, \emph{On varieties of minimal degree (a centennial account)}, Proc. of Symposia in Pure Math. \textbf{46}, 3--43 (1987)
\bibitem{[GS06]}
M. Gross and B. Siebert, \emph{Mirror symmetry via logarithmic degeneration data I}, J. Diff. Geom \textbf{72}, 169--338 (2006)
\bibitem{[GS10]}
M. Gross and B. Siebert, \emph{Mirror symmetry via logarithmic degeneration data II}, J. Alg. Geom. \textbf{19}, 679--780 (2010)
\bibitem{[JK96]}
T. Johnsen and S. Kleiman, \emph{Rational curves of degree at most 9 on a general quintic threefold}, Comm. Algebra. \textbf{24} (8), 2721--2753 (1996)
\bibitem{[Ka86a]}
S. Katz, \emph{On the finiteness of rational curves on quintic threefolds}, Compos. Math. \textbf{60}, 151--162 (1986)
\bibitem{[Ka86b]}
S. Katz, \emph{Lines on Complete Intersection Threefolds with K=0}, Math. Z. \textbf{191}, 293--296 (1986)
\bibitem{[Kn12]}
A. L. Knutsen, \emph{On isolated smooth curves of low genera in Calabi-Yau complete intersection threefolds}, Trans. Amer. Math. Soc. \textbf{364}, 5243--5264 (2012)
\bibitem{[Li01]}
J. Li, \emph{Stable morphisms to singular schemes and relative stable morphisms}, J. Diff. Geom. \textbf{57}, 509--578 (2001)
\bibitem{[Li02]}
J. Li, \emph{A degeneration formula of GW-invariants}, J. Diff. Geom. \textbf{60}, 199-293 (2002)
\bibitem{[Ni15]}
T. Nishinou, \emph{Counting curves via degenerations}, Preprint, \url{http://arxiv.org/pdf/1108.2802v3.pdf} (2015)
\bibitem{[Ni95]}
P. Nijsse, \emph{Clemens conjecture for octic and nonic curves}, Indag. Math. \textbf{6} (2), 213--221 (1995)

\end{thebibliography}
\end{document}